\documentclass[preprints,article,accept,moreauthors,pdftex]{my-mdpi}

\firstpage{1}
\pubvolume{11}
\issuenum{6}
\articlenumber{295}
\doinum{10.3390/axioms11060295}
\pubyear{2022}
\copyrightyear{2022}
\externaleditor{Academic Editor: Chris Goodrich}
\datereceived{2 May 2022}
\daterevised{7 June 2022} 
\dateaccepted{14 June 2022}
\datepublished{16 June 2022}
\hreflink{https://doi.org/\\10.3390/axioms11060295}

\pdfoutput=1

\Title{Existence Results for a Multipoint Fractional Boundary Value Problem in the Fractional Derivative Banach Space}

\TitleCitation{Existence Results for a Multipoint Fractional Boundary Value Problem in the Fractional Derivative Banach Space}

\Author{Djalal Boucenna $^{1,\dagger}$\orcidA{}, 
Amar Chidouh $^{2,\dagger}$\orcidB{} 
and Delfim F. M. Torres $^{3,}$*$^{,\dagger}$\orcidC{}}

\AuthorNames{Djalal Boucenna, Amar Chidouh and Delfim F. M. Torres}

\AuthorCitation{Boucenna,~D.; Chidouh,~A.; Torres,~D.F.M.}

\address{$^{1}$ \quad High School of Technological Teaching, Enset, Skikda 21001, Algeria; 
djalal.boucenna@enset-skikda.dz\\
$^{2}$ \quad Department of Mathematics, Chadli Bendjedid University, Eltarf 36000, Algeria; 
m2ma.chidouh@gmail.com\\
$^{3}$ \quad Center for Research and Development in Mathematics and Applications (CIDMA), 
\mbox{{Department of Mathematics,}} University of Aveiro, 3810-193 Aveiro, Portugal}

\corres{Correspondence: delfim@ua.pt}

\firstnote{These authors contributed equally to this work.}

\abstract{{We study a class of nonlinear implicit fractional differential 
equations subject to nonlocal boundary conditions expressed in terms 
of nonlinear integro-differential equations. Using the Krasnosel'skii fixed-point theorem  
we prove, via the Kolmogorov--Riesz criteria, the existence of solutions. The existence results 
are established in a specific fractional derivative Banach space and they are 
illustrated by two numerical examples.}}



\begin{document}


\section{Introduction}

It is noticeable, in recent years, that the field of fractional calculus
has been swept for research by many  mathematicians, due to its
effectiveness in describing many physical phenomena, see,
e.g., {\cite{rev1,rev2,rev3,rev4,rev5,rev6,MR4200529}}.

A fractional derivative is a generalization of the ordinary one.
Despite the emergence of several definitions of fractional derivative,
the content is one that depends entirely on Volterra integral equations
and their kernel, which facilitates the description of each phenomenon
 {as} a temporal lag, such as rheological phenomena \cite{12,13,11}.

The study of differential equations is considered of primary importance
in mathematics. In applications, differential equations serve as
mathematical models for all natural phenomena.
Regardless of their type (ordinary, partial, or fractional),
the study of differential equations is developed in three directions:
existence, uniqueness, and stability of solutions.  Therefore, to
investigate boundary value problems is always a central question in mathematics
 {\cite{rev7,rev8,rev9,MR4014986,MR3971207,MR4221136}}.

Often, it is of central importance to know the behavior of any solution,
of the equation under study, at the boundary  of the domain, because
that makes it easier to find the solution. In 2009, Ahmad and Nieto
considered the following boundary value problem \cite{14}:
\begin{equation}
\label{66}
\begin{gathered}
^{C}D^{\alpha}y\left(  t\right)  =f\left(  t,y\left(  t\right),
\int_{0}^{t}\varphi(t,s)y(s)ds\right), \quad 1<\alpha<2,\\
a\text{ }y\left(  0\right)  +by^{\prime}\left(  0\right)
=\int_{0}^{1} q_{1}(y\left(  s\right)  )ds\text{,}\\
a y\left(  1\right)  +by^{\prime}\left(  1\right)
=\int_{0}^{1}q_{2}(y\left(s\right)  )ds\text{.}
\end{gathered}
\end{equation}
Their results of existence are obtained via Krasnosel'skii fixed-point
theorem in the space of continuous functions. For that, they apply Ascoli's
theorem in order to provide the compactness of the first part of the
Krasnosel'skii operator.

The pioneering work of Ahmad and Nieto of 2009 \cite{14}
gave rise to several different investigations. These include:
inverse source problems for fractional integrodifferential equations \cite{MR3272530};
the study of positive solutions for singular fractional boundary value problems
with coupled integral boundary conditions \cite{MR3278339};
the expression and properties of Green's function for nonlinear boundary value problems
of fractional order with impulsive differential equations \cite{MR3481182};
existence of solutions to several kinds of differential equations using
the coincidence theory \cite{MR3434271}; existence and uniqueness of solution
for fractional differential equations with Riemann--Liouville fractional integral
boundary conditions \cite{MR3442607}; sufficient conditions for the existence
and uniqueness of solutions for a class of terminal value problems
of fractional order on an infinite interval \cite{MR3416395}; existence of solutions,
and stability, for fractional integro-differential equations involving a general form
of Hilfer fractional derivative with respect to another function \cite{MR3911898};
existence of solutions for a boundary value problem involving mixed generalized
fractional derivatives of Riemann--Liouville and Caputo, supplemented
with nonlocal multipoint boundary conditions \cite{MR4178811};
existence conditions to fractional order hybrid differential equations \cite{MR4185337};
and an existence analysis for a nonlinear implicit fractional differential equation
with integral boundary conditions \cite{MR4201632}.
Motivated by all these existence results, we consider here
a more general multipoint fractional boundary value problem
in the fractional derivative Banach space.

Let $1<p<\infty$ and $1\geq\gamma>\frac{1}{p}$
and consider the following fractional boundary \mbox{value problem:}
\begin{equation}
\label{77}
\begin{gathered}
^{C}D^{\alpha}y\left(  t\right)  =f\left(  t,y\left(  t\right),
^{C}D^{\gamma}y\left(  t\right)  \right)
+^{C}D^{\alpha-2}g\left(  t,y\left(  t\right)\right), \quad 2<\alpha<3,\\
\text{ }y\left(  0\right)  +y^{\prime}\left(  0\right)
=\int_{0}^{1} q_{1}(y(s))ds,\\
y\left(  1\right)  +y^{\prime}\left(  1\right)
=\int_{0}^{1}q_{2}(y(s))ds,\\
y^{\prime\prime}\left( 0\right)=0,
\end{gathered}
\end{equation}
where $^{C}D^{\alpha}$ is the standard Caputo derivative,
$f:\left[0,1\right]  \times\mathbb{R}\times\mathbb{R}\rightarrow\mathbb{R}$,
and $g:\left[  0,1\right]  \times\mathbb{R}\rightarrow\mathbb{R}$
and $q_{1},q_{2}:\mathbb{R}\rightarrow\mathbb{R}$ are given functions such that
$g(t,0) {=g(0,y)}=q_{1}(0)=q_{2}(0)=0$ for any $(t,y) \in\left[0,1\right]  \times\mathbb{R}$.
 {Our problem \eqref{77} generalizes \eqref{66}
and finds applications in viscoelasticity, where the 
fractional operators are associated with delay kernels that make the fractional 
differential equations the best models for several rheological Maxwell phenomena.
In particular, for $\alpha \in (1,2)$, we can model oscillatory processes with 
fractional damping \cite{Bagley:Torvik}.}

We prove existence of a solution to problem \eqref{77}
in the special Banach space $E^{\gamma,p}$ {that is known in the literature
as the fractional derivative space \cite{15}. This Banach space is equipped with the norm}
\begin{equation}
\label{n}
\left\Vert u\right\Vert _{\gamma,p}
=\left(  \int_{0}^{T}\left\vert u\left(t\right)  \right\vert ^{p}
+\int_{0}^{T}\left\vert ^{C}D_{0}^{\gamma}u\left(t\right)
\right\vert ^{p}\right)^{\frac{1}{p}}{.}
\end{equation}

The paper is organized as follows. In Section~\ref{sec2},
we recall some useful definitions and lemmas to prove our main results.
The original contributions are then given in Section~\ref{sec3}. The main
result is Theorem~\ref{the3.1}, which establishes the existence of solutions
to the fractional boundary value problem \eqref{77} using Krasnosel'skii
fixed point theorem.  {Two illustrative examples are given.} We end with
Section~\ref{sec:disc}, discussing the obtained existence result.


\section{Preliminaries}
\label{sec2}

For the convenience of the reader, and to facilitate the analysis of
problem \eqref{77}, we begin by recalling the necessary background
from the theory of fractional calculus \cite{kilbas,podlubny}.

\begin{Definition}
\label{def1}
The Riemann--Liouville fractional integral of order $\alpha>0$ of
a function \ $f:\left(  0,+\infty\right)  \rightarrow \mathbb{R}$
is given by
\[
I_{0}^{\alpha}f\left(  t\right)  =\frac{1}{\Gamma\left(  \alpha\right)  }
\int_{0}^{t}\left(  t-s\right)  ^{\alpha-1}f\left(  s\right)  ds{.}
\]
\end{Definition}

\begin{Definition}
\label{def2}
The Caputo fractional derivative of order $\alpha>0$ of a
function $f:\left(  0,+\infty\right)  \rightarrow \mathbb{R}$
is given by
\[
^{C}D_{0}^{\alpha}f\left(  t\right)  =\frac{1}{\Gamma\left(  n-\alpha\right)
}\int_{0}^{t}\frac{f^{\left(  n\right)  }\left(  s\right)  }{\left(
t-s\right)  ^{\alpha-n+1}}ds=I_{0}^{n-\alpha}f^{\left(  n\right)  }\left(
t\right),
\]
where $n=\left[  \alpha\right]  +1$, with $\left[  \alpha\right]$
denoting the integer part of $\alpha$.
\end{Definition}

\begin{Lemma}[See \cite{14}]
\label{lem1}
Let $\alpha>0$. Then, the fractional differential equation 
 {$^{C}D_{0^{+}}^{\alpha}u\left(  t\right) =0$ has
\[
u\left(  t\right)  =c_{0}+c_{1}t+c_{2}t^{2}+\cdots+c_{n-1}t^{n-1},
\quad c_{i}\in \mathbb{R}, \quad i=1,2,\ldots,n-1,
\]
as solution.}
\end{Lemma}

\begin{Definition}
A map $f:\left[  0,1\right]
\times\mathbb{R}\times\mathbb{R}\rightarrow \mathbb{R}$
is said to be Carath\'{e}odory if
\begin{enumerate}
\item[(a)] $t\rightarrow f\left(  t,u;v\right)$
is measurable for each $u,v\in \mathbb{R}$;

\item[(b)] $\left(  u,v\right)  \rightarrow f\left(  t,u;v\right)$
is continuous for almost all $t\in\left[  0,1\right]$.
\end{enumerate}
\end{Definition}

\begin{Definition}
\label{def4}
Let $J$ be a measurable subset of $\mathbb{R}$
and $f:J\times \mathbb{R}^{d_{1}}\rightarrow\mathbb{R}^{d_{2}}$
satisfies the Carath\'{e}odory condition. By a generalized
Nemytskii operator we mean the mapping $N_{f}$ taking a (measurable) vector
function $u=\left(  u_{1},\ldots,u_{d_{1}}\right)  $ to the function
$N_{f}u(t)=f(t,u(t))$, $t\in J$.
\end{Definition}

The following lemma is concerned with the continuity
of the operator $N_{f}$ with $d_{1}=2$ and $d_{2}=1$.

\begin{Lemma}[See \cite{precup}]
\label{lem2}
Consider the same data of Definition~\ref{def4}. Assume there exists
$w\in L^{1}\left(  [0,1]\right)  $ with $1\leq p<\infty$ and a constant $c>0$
such that
$\left\vert f\left(  t,u,v\right)  \right\vert \leq w\left(  t\right)
+c\left(  \left\vert u\right\vert ^{p}+\left\vert v\right\vert ^{p}\right)$
for almost all $t\in\lbrack0,1]$ and $u,v\in \mathbb{R}$.
Then, the Nemytskii operator
$$
N_{f}u(t)=f\left(  t,u(t)\right),
\quad
u=\left(  u_{1},u_{2}\right)  \in L^{p}\left(  0,1\right)
\times L^{p}\left(  0,1\right),
\quad t\in\lbrack0,1] \text{ a.e.},
$$
is continuous from $L^{p}\left([0,1]\right)  \times L^{p}\left(  [0,1]\right)$
to $L^{1}\left(  0,1\right)$.
\end{Lemma}

\begin{Lemma}[See \cite{brizis}]
\label{lem3} Let $\mathcal{F}$
be a bounded set in $L^{p}\left(  [0,1]\right)$
with $1\leq p<\infty$. Assume that
\[
\underset{\left\vert h\right\vert \rightarrow0}{\lim}\left\Vert \tau
_{h}f-f\right\Vert _{p}=0\ \text{uniformly\ on}\ \mathcal{F}.
\]
Then, $\mathcal{F}$ is relatively compact in $L^{p}\left(  [0,1]\right)$.
\end{Lemma}

For any $1\leq p<\infty$, we denote
\[
\left\Vert u\right\Vert _{L^{p}\left[  0,T\right]  }:=\left(
\int_{0}^{T}\left\vert u\left(  t\right)  \right\vert ^{p}\right)^{\frac{1}{p}},
\quad  \left\Vert u\right\Vert _{\infty}
:=\max_{t\in\left[  0,T\right]  }\left\vert u\left(  t\right)  \right\vert.
\]
Now, we give the definition and some properties of $E^{\gamma,p}$. For
more details about the following lemmas, see \cite{15,nyamoradi} and
references therein.

\begin{Definition}
\label{def3}
Let $0<\gamma\leq1$ and $1<p<\infty$. The fractional derivative space
$E^{\gamma,p}$ is defined by the closure of $C^{\infty}\left(
\left[  0,T\right]\right)$ with respect to the norm
\begin{equation}
\label{2.1}
\left\Vert u\right\Vert _{\gamma,p}=\left(  \int_{0}^{T}\left\vert u\left(
t\right)  \right\vert ^{p}+\int_{0}^{T}\left\vert ^{C}D_{0}^{\gamma}u\left(
t\right)  \right\vert ^{p}\right)  ^{\frac{1}{p}}.
\end{equation}
\end{Definition}

\begin{Lemma}[See \cite{15,nyamoradi}]
\label{lem4}
Let $0<\gamma\leq1$ and $1<p<\infty$. The fractional derivative
space $E^{\gamma,p}$ is a reflexive and separable Banach space.
\end{Lemma}

\begin{Lemma}[See \cite{15,nyamoradi}]
\label{lem5}
Let $0<\gamma\leq1$ and $1<p<\infty$. For all $u\in E^{\gamma,p}$, we have
\begin{equation}
\label{2.2}
\left\Vert u\right\Vert _{L^{p}}\leq\frac{T^{\alpha}}{\Gamma\left(
\gamma+1\right)  }\left\Vert ^{C}D_{0}^{\gamma}u\right\Vert _{L^{p}}.
\end{equation}
Moreover, if $\gamma>\frac{1}{p}$ and $\frac{1}{p}+\frac{1}{q}=1$, then
\begin{equation}
\label{2.3}
\left\Vert u\right\Vert _{\infty}\leq\frac{T^{\alpha-\frac{1}{p}}}{\Gamma\left(
\gamma\right)  \left(  \left(  \gamma-1\right)  q+1\right)^{\frac{1}{q}}}
\left\Vert ^{C}D_{0}^{\gamma}u\right\Vert _{L^{p}}.
\end{equation}
\end{Lemma}

According to the inequality \eqref{2.2}, we can also consider the space
$E^{\gamma,p}$ with respect to the equivalent norm
\[
\left\Vert u\right\Vert _{\gamma,p}=\left\Vert ^{C}D_{0}^{\gamma}u\right\Vert
_{L^{p}}=\left(  \int_{0}^{T}\left\vert ^{C}D_{0}^{\gamma}u\left(  t\right)
\right\vert ^{p}\right)  ^{\frac{1}{p}}{,\,\,\, u\in E^{\gamma,p}.}
\]


\section{Main Results}
\label{sec3}

We begin by considering a linear problem and obtain its solution
in terms of a \mbox{Green function.}

\begin{Lemma}
Assume $h,k\in C\left( \left[  0,1\right]\right)$,
${k(0)=0}$ and $\alpha \in (2,3)$.
Then, the solution to the boundary value problem
\begin{equation}
\label{2.4}
\begin{gathered}
^{C}D^{\alpha}y\left(  t\right)  =h\left(  t\right)
+^{C}D^{\alpha-2}k\left(t\right)  \text{,}
\quad t\in\left(  0,1\right)  \text{,}\\
y^{\prime\prime}\left(  0\right)  =0\text{,}\\
\text{ }y\left(  0\right)  +y^{\prime}\left(  0\right)  =\int_{0}^{1}\eta
_{1}\left(  s\right)  ds\text{,}\\
y\left(  1\right)  +y^{\prime}\left(  1\right)  =\int_{0}^{1}\eta_{2}\left(
s\right)  ds\text{,}
\end{gathered}
\end{equation}
is given by
\begin{equation*}
y\left(  t\right)   =\int\limits_{0}^{1}G\left(  t,s\right)  h\left(
s\right)  ds+  \int_{0}^{1}H\left(  t,s\right) k\left(s\right)  ds
+\left(  2-t\right)  \int_{0}^{1}\eta_{1}\left(  s\right)  ds+\left(
t-1\right)  \int_{0}^{1}\eta_{2}\left(  s\right)  ds\text{,}
\end{equation*}
where
\begin{equation}
\label{eq:G}
G\left(  t,s\right)  =\left\{
\begin{array}[c]{ll}
\frac{\left(  t-s\right)  ^{\alpha-1}+\left(  1-t\right)  \left(  1-s\right)
^{\alpha-1}}{\Gamma\left(  \alpha\right)  }+\frac{\left(  1-t\right)  \left(
1-s\right)  ^{\alpha-2}}{\Gamma\left(  \alpha-1\right)  }{,} & 0\leq s\leq
t\leq1\text{,}\vspace{0.3cm}\\
\frac{\left(  1-t\right)  \left(  1-s\right)  ^{\alpha-1}}{\Gamma\left(
\alpha\right)  }+\frac{\left(
1-t\right)\left(  1-s\right)^{\alpha-2}}{\Gamma\left(  \alpha-1\right)  }{,}
& 0\leq t\leq s\leq1{,}
\end{array}
\right.
\end{equation}
and
\begin{equation}
\label{eq:H}
H\left(  t,s\right)  =\left\{
\begin{array}[c]{ll}
\left(t-s \right) + \left(  1-t\right)\left(  2-s\right){,} & 0\leq s\leq
t\leq1\text{,}\vspace{0.3cm}\\
\left(  1-t\right)\left(  2-s\right){,}
& 0\leq t\leq s\leq1\text{.}
\end{array}
\right.
\end{equation}
\end{Lemma}

\begin{proof}
Let $y$ be a solution of problem \eqref{2.4}. By Lemma~\ref{lem1},
we have
\[
y\left(  t\right)  =c_{0}+c_{1}t+c_{2}t^{2}+\frac{1}{\Gamma\left(
\alpha\right)  }\int_{0}^{t}\left(  t-s\right)  ^{\alpha-1}h\left(  s\right)
ds+I_{0}^{2}k\left(  t\right).
\]
Taking the conditions \eqref{2.4} into account, it follows that
\[
c_{2}=0\text{,}
\]
\[
y\left(  0\right)  +y^{\prime}\left(  0\right)  =c_{0}+c_{1}
=\int_{0}^{1} \eta_{1}\left(  s\right)  ds\text{,}
\]
and
\begin{align*}
y\left(  1\right)  +y^{\prime}\left(  1\right)
&=c_{0}+2c_{1}+\frac{1}{\Gamma\left(  \alpha\right)}
\int_{0}^{1}\left(  1-s\right)  ^{\alpha-1}h\left(  s\right)  ds
+\int_{0}^{1}\left(  1-s\right)  k\left(  s\right) ds\\
& \quad +\frac{1}{\Gamma\left(  \alpha-1\right)  }
\int_{0}^{1}\left(  1-s\right)^{\alpha-2}h\left(  s\right)  ds
+\int_{0}^{1}k\left(  s\right)  ds\\
&  =\int_{0}^{1}\eta_{2}\left(  s\right)  ds,
\end{align*}
which implies
\begin{align*}
c_{0}  &  =\frac{1}{\Gamma\left(  \alpha\right)  }\int_{0}^{1}\left(
1-s\right)  ^{\alpha-1}h\left(  s\right)  ds+\frac{1}{\Gamma\left(
\alpha-1\right)  }\int_{0}^{1}\left(  1-s\right)  ^{\alpha-2}h\left(
s\right)  ds\\
& \quad +\int_{0}^{1}\left(  2-s\right)  k\left(  s\right)  ds+2\int_{0}^{1}
\eta_{1}\left(  s\right)  ds-\int_{0}^{1}\eta_{2}\left(  s\right)  ds{,}
\end{align*}
and
\begin{align*}
c_{1}  &  =-\frac{1}{\Gamma\left(  \alpha\right)  }\int_{0}^{1}\left(
1-s\right)  ^{\alpha-1}h\left(  s\right)  ds-\frac{1}{\Gamma\left(
\alpha-1\right)  }\int_{0}^{1}\left(  1-s\right)  ^{\alpha-2}h\left(
s\right)  ds\\
&  \quad -\int_{0}^{1}\left(  2-s\right)  k\left(  s\right)  ds
+\int_{0}^{1}\eta_{2}\left(  s\right)  ds
-\int_{0}^{1}\eta_{1}\left(  s\right)  ds\text{.}
\end{align*}
Hence, the solution of problem \eqref{2.4} is
\begin{align*}
y\left(  t\right)
&=\int_{0}^{t}\left(\frac{\left(
t-s\right)^{\alpha-1}+\left(  1-t\right)
\left(  1-s\right)^{\alpha-1}}{\Gamma\left(
\alpha\right)  }+\frac{\left(  1-t\right)
\left(  1-s\right)^{\alpha-2}}{\Gamma\left(  \alpha-1\right)  }\right)
h\left(  s\right)  ds\\
&\quad  +\int_{t}^{1}\left(  \frac{\left(  1-t\right)
\left(  1-s\right)^{\alpha-1}}{\Gamma\left(  \alpha\right)  }
+\frac{\left(  1-t\right)  \left(
1-s\right)  ^{\alpha-2}}{\Gamma\left(  \alpha-1\right)  }\right)  h\left(
s\right)  ds\\
& \quad +  \int_{0}^{t}\left(  \left(  t-s\right)
+\left(  1-t\right)\left(  2-s\right)\right)   k (s)ds
+ \int_{t}^{1}\left(  1-t\right)\left(  2-s\right)  k\left(  s\right)ds\\
&  \quad + \left(  2-t\right)  \int_{0}^{1}\eta_{1}\left(  s\right)  ds
+\left(  t-1\right)  \int_{0}^{1}\eta_{2}\left(  s\right)  ds\\
& =\int_{0}^{1}G\left(  t,s\right)  h\left(  s\right)  ds+
\int_{0}^{1} H\left(  t,s\right) k\left(  s\right)  ds\\
& \quad +\left(  2-t\right)  \int_{0}^{1}\eta_{1}\left(  s\right)  ds+\left(
t-1\right)  \int_{0}^{1}\eta_{2}\left(  s\right)  ds\text{.}
\end{align*}
The proof is complete.
\end{proof}

\begin{Lemma}
Functions $G, H , \frac{\partial^{\gamma}}{\partial t}G $
and $\frac{\partial^{\gamma}}{\partial t}H$
are continuous on $\left[  0,1\right]
\times\left[  0,1\right]$ and satisfy  {for all $t,s\in [0,1]$:}
\begin{enumerate}
\item
$
\left\vert G\left(  t,s\right)  \right\vert
\leq\frac{3}{\Gamma\left(  \alpha-1\right)},
\quad \left\vert H\left(  t,s\right)  \right\vert \leq 3{.}
$
\item
$
\left\vert \frac{\partial^{\gamma}}{\partial t}
G\left(  t,s\right)  \right\vert
\leq\frac{\Gamma\left(\alpha \right) }{\Gamma\left(
\alpha-\gamma\right)  }+\frac{2}{\Gamma(2-\gamma)\Gamma\left(
\alpha-1\right)  },  \quad \left\vert \frac{\partial^{\gamma}}{\partial t}
H\left(  t,s\right)  \right\vert \leq\frac{3}{\Gamma(2-\gamma)}\allowbreak
.$
\end{enumerate}
\end{Lemma}

\begin{proof}
We have
\begin{equation}
\label{a10}
^{C}D_{0}^{\gamma}(1-t)= I_{0}^{1-\gamma}(1-t)^{\prime}
=  -\frac{1}{\Gamma({2-\gamma})}t^{1-\gamma}{,}
\end{equation}
and
\[
\frac{\partial^{\gamma}}{\partial t} (t-s)^{\alpha-1} = I_{0}^{1-\gamma}
\frac{\partial}{\partial t} (t-s)^{\alpha-1} = (\alpha-1) I_{0}^{1-\gamma}
(t-s)^{\alpha-2}.
\]
Thus, for $0\leq s\leq t\leq1$, we get
$\frac{\partial^{\gamma}}{\partial t}(t-s)^{\alpha-1} \geq 0$ and
\begin{equation}
\label{a11}
\frac{\partial^{\gamma}}{\partial t} (t-s)^{\alpha-1}
\leq ^{C}D_{0}^{\gamma} t^{\alpha-1}
= \frac{\Gamma\left(\alpha \right)}{\Gamma\left(
\alpha- \gamma\right)}t ^{\alpha- \gamma-1}.
\end{equation}
On the other hand, we have $\Gamma\left(  \alpha-1 \right)
\leq\Gamma\left(\alpha\right)$ for $2 \leq\alpha\leq 3$.
Now, we give the bound of functions $\left|  G\left(  t,s\right)  \right|$
and $\left|  \frac {\partial^{\gamma}}{\partial t}G\left(  t,s\right)  \right|$.
From the definition of function $G$ and \eqref{a10} and \eqref{a11}, \mbox{we obtain:}
\begin{itemize}
\item For $0\leq s\leq t\leq1$,
\[
\begin{array}[c]{lll}
\left|  G\left(  t,s\right)  \right|  & = & \frac{\left(  t-s\right)
^{\alpha-1}+\left(  1-t\right)  \left(  1-s\right)  ^{\alpha-1}}{\Gamma\left(
\alpha\right)  }+\frac{\left(  1-t\right)  \left(  1-s\right)^{\alpha-2}}{
\Gamma\left(  \alpha-1\right)  }\\
& \leq & \frac{\left(  1-s\right)  ^{\alpha-1}\left(  1 +\left(  1-t\right)
\right)  }{\Gamma\left(  \alpha\right)  }+\frac{\left(  1-t\right)  \left(
1-s\right)  ^{\alpha-2}}{\Gamma\left(  \alpha-1\right)  }\\
& \leq & \frac{ 1 +\left(  1-t\right)  }{\Gamma\left(  \alpha\right)}
+\frac{\left(  1-t\right)  }{\Gamma\left(  \alpha-1\right)}\\
& \leq & \frac{3}{\Gamma\left(  \alpha-1 \right)}{,}
\end{array}
\]
and
\[
\begin{array}[c]{lll}
\left|  \frac{\partial^{\gamma}}{\partial t}G\left(  t,s\right)  \right|
&\leq & \left|  \frac{\Gamma\left(  \alpha\right) }{\Gamma\left(  \alpha
-\gamma\right)  }t^{\alpha- \gamma-1} \right|
+ \left|  \frac{t^{1-\gamma} }{\Gamma(2-\gamma)}\left(
\frac{ (1-s) ^{\alpha-1}}{\Gamma( \alpha)}+\frac{ \left(
1-s\right)  ^{\alpha-2}}{\Gamma\left(  \alpha-1\right)  }\right)  \right| \\
& \leq & \frac{{\Gamma\left(\alpha \right) }}{\Gamma\left(  \alpha- \gamma\right)  }
+\frac{1}{\Gamma(2-\gamma)}\left(  \frac{ 1}{\Gamma( \alpha)}
+\frac{ 1}{\Gamma\left(\alpha-1\right)  }\right) \\
& \leq & \frac{{\Gamma\left(\alpha \right) }}{\Gamma\left(\alpha- \gamma\right)}
+\frac{2}{\Gamma(2-\gamma)\Gamma\left(  \alpha-1\right)}.
\end{array}
\]

\item For $0\leq t\leq s\leq1$,
\[
\begin{array}[c]{lll}
\left|  G\left(  t,s\right)  \right|  & = & \frac{\left(  1-t\right)  \left(
1-s\right)  ^{\alpha-1}}{\Gamma\left(  \alpha\right)  }+\frac{\left(
1-t\right)  \left(  1-s\right)  ^{\alpha-2}}{\Gamma\left(  \alpha-1\right)
}\\
& \leq & \frac{ \left(  1-t\right)  }{\Gamma\left(  \alpha\right)}
+\frac{\left(  1-t\right)  }{\Gamma\left(  \alpha-1\right)  }\\
& \leq & \frac{2 }{\Gamma\left(  \alpha-1 \right)}{,}
\end{array}
\]
and
\[
\begin{array}[c]{lll}
\left|  \frac{\partial^{\gamma}}{\partial t}G\left(  t,s\right)  \right|
&=& \left|  -\frac{t^{1-\gamma} }{\Gamma(2-\gamma)}\left(
\frac{ (1-s)^{\alpha-1}}{\Gamma( \alpha)}+\frac{ \left(  1-s\right)^{\alpha-2}}{
\Gamma\left(  \alpha-1\right)  }\right)  \right| \\
& \leq & \frac{1}{\Gamma(2-\gamma)}\left(  \frac{ 1}{\Gamma(
\alpha)}+\frac{ 1}{\Gamma\left(  \alpha-1\right)  }\right) \\
& \leq & \frac{2 }{\Gamma(2-\gamma)\Gamma\left(  \alpha-1\right)}.
\end{array}
\]
\end{itemize}
By using the same above calculation, we obtain
the estimation of $\left| H(t,s)\right|$
and $\left| \frac{\partial^{\gamma}}{\partial t}
H(t,s)\right|$. The proof is complete.
\end{proof}

In the sequel, we denote
$$
G_{\gamma}\left(  t,s\right)
:=\frac{\partial^{\gamma}}{\partial t} G\left(  t,s\right),
\quad t,s\in\left[  0,1\right]\times\left[  0,1\right].
$$
Moreover, we also use the following notations:
$G^{\ast}:=\max_{t,s\in\left[  0,1\right]
\times\left[  0,1\right]  }\left|  G\left(  t,s\right)\right|$ and
$$
G_{\gamma}^{\ast}:=\max_{t,s\in\left[  0,1\right] 
\times\left[  0,1\right]}\left|
G_{\gamma} \left(  t,s\right)  \right|.
$$

\begin{Theorem}
\label{the3.1}
Assume that the following four hypotheses hold:
\begin{enumerate}

\item[(H1)] $f:\left[  0,1\right]  \times\mathbb{R}
\times\mathbb{R}\rightarrow\mathbb{R}$
satisfies the Carath\'{e}odory condition.

\item[(H2)] There exist $w\in L^{1}\left(  0,1\right)$ and $c>0$ such that
\begin{equation}
\label{3.1}
\left\vert f\left(  t,u,v\right)  \right\vert \leq w\left(  t\right)
+c\left(  \left\vert u\right\vert ^{p}+\left\vert v\right\vert ^{p}\right)
\text{ for }t\in\left(  0,1\right)  \text{ and }u,v\in\mathbb{R}\text{.}
\end{equation}

\item[(H3)] {There} exist two strictly positive constants $k_1$ and $k_2$
and a function $\varphi_{1} \in L^{q}\left(  \left(  0,1\right),
\mathbb{R}_{+}\right)$, $\frac{1}{p}+\frac{1}{q}=1$,
such that {for all $t\in\left[0,1\right]$ and $x,y\in\mathbb{R}$, we have}
\begin{equation*}
\begin{split}
\left\vert g\left(  t,x\right)  -g\left(  t,y\right)  \right\vert
&\leq \varphi_{1}\left(  t\right)  \left\vert x-y\right\vert,\\
\left\vert q_{1}\left(  x\right)  -q_{1}\left(  y\right)  \right\vert
&\leq k_{1} \left\vert x-y\right\vert,\\
\left\vert q_{2}\left(  x\right)  -q_{2}\left(  y\right)  \right\vert
&\leq k_{2}\left\vert x-y\right\vert{.}
\end{split}
\end{equation*}

\item[(H4)] There exists a real number $R>0$ such that
\begin{equation}
\label{3.2}
\frac{R\left[  3\left\Vert \varphi_{1}\right\Vert _{q}
+k_{1}+k_{2}\right]}{
\Gamma\left(  2-\gamma\right)  \Gamma\left(  1+\gamma\right)  }
+G_{\gamma}^{\ast}\left(\left\Vert w\right\Vert_{1}+c\left(
1+\left( \frac{1}{\Gamma( \gamma+1)}\right)^{p}\right)
R^{p}\right)  \leq R\text{.}
\end{equation}
\end{enumerate}
{Then, if}
\begin{equation}
\label{3.3}
\frac{ 3\left\Vert \varphi_{1}\right\Vert _{q}+ k_1
+ k_2 }{\Gamma\left(2-\gamma\right)
\Gamma\left(  1+\gamma\right)  }<1,
\end{equation}
the boundary value problem \eqref{77} has a solution in $E^{\gamma,p}$.
\end{Theorem}

\begin{proof}
We transform problem \eqref{77} into a fixed-point problem. Define two
operators $F,L :E^{\gamma,p}\rightarrow E^{\gamma,p}$ by
\[
Fy\left(  t\right)  =\int_{0}^{1}G\left(  t,s\right)
f\left(  s,y\left(s\right),
D^{\gamma}y\left(  s\right)  \right)  ds{,}
\]
and
\begin{equation*}
Ly\left(  t\right)   =  \int_{0}^{1}H\left(  t,s\right)
g\left(  s,y\left(  s\right)  \right)  ds+\left(  2-t\right)
\int_{0}^{1}q_{1}\left(  y\left(  s\right)  \right)  ds
+\left(  t-1\right)  \int_{0}^{1}q_{2}\left(  y\left(  s\right)  \right)
ds\text{.}
\end{equation*}
Then, $y$ is a solution of problem \eqref{77} if, and only if, $y$ is a fixed
point of $F+L$. We define the set $B_{R}$ as follows:
\[
B_{R}:=\left\{  u\in E^{\gamma,p},\left\Vert
u\right\Vert_{E^{\gamma,p}}\leq R\right\}  \text{,}
\]
where $R$ is the same constant defined in $\left(  H_{3}\right)$. It is
clear that $B_{R}$ is convex, closed, and a bounded subset of $E^{\gamma,p}$.
We shall show that $F,G$ satisfy the assumptions of Krasnosel'skii 
fixed-point theorem. The proof is given in several steps.
\begin{enumerate}
\item[(i)] We prove that $F$ is continuous. Let $\left(  y_{n}\right)_{n\in\mathbb{N}}$
be a sequence such that $y_{n}\rightarrow y$ in $E^{\gamma,p}$. From \eqref{3.1}
and Lemma~\ref{lem2}, and for each $t\in\lbrack0,1]$, we obtain
\begin{align*}
& \left\vert \left(  ^{C}D_{0}^{\gamma}Fy_{n}\right)  \left(  t\right)
-\left(  ^{C}D_{0}^{\gamma}Fy\right)  \left(  t\right)  \right\vert \\
& \leq \int_{0}^{1}\left\vert G_{\gamma}\left(  t,s\right)  \right\vert
\,
\left\vert f\left(  s,y_{n}\left(  s\right)  ,D^{\gamma}y_{n}\left(  s\right)
\right)  -f\left(  s,y\left(  s\right)  ,D^{\gamma}y\left(  s\right)  \right)
\right\vert ds \\
& \leq G_{\gamma}^{\ast}\left\Vert N_{f}y_{n}-N_{f}y\right\Vert_{1}\text{.}
\end{align*}
Applying the $L^{p}$ norm, we obtain that
$\left\Vert Fy_{n}-Fy\right\Vert _{E^{\gamma,p}}\rightarrow 0$
when $y_{n}\rightarrow y\ \text{in}E^{\gamma,p}$.
Thus, the operator $F$ is continuous.

\item[(ii)] Now, we prove that $Fx+Ly\in B_{R}$ for $x,y\in B_{R}$.
Let $x,y\in B_{R}$, $t\in\left(0,1\right)$. In view of
hypothesis $\left(H3\right)$, we obtain
\begin{align*}
\left\vert ^{C}D_{0}^{\gamma}Fy\left(  t\right)  \right\vert
&  \leq\int_{0}^{1}\left\vert G_{\gamma}\left(  t,s\right)  \right\vert
\left\vert f\left(  s,y\left(  s\right)  ,D^{\gamma}y\left(  s\right)  \right)
\right\vert ds\\
&  \leq G_{\gamma}^{\ast}\left(  \left\Vert w\right\Vert _{1}+c\left(
\left\Vert y\right\Vert _{p}^{p}+\left\Vert ^{C}D_{0}^{\gamma}y\right\Vert
_{p}^{p}\right)  \right)  \\
&  \leq G_{\gamma}^{\ast}\left(  \left\Vert w\right\Vert _{1}+c\left(
1+\left( \frac{1}{\Gamma(\gamma+1)}\right) ^{p}\right)  R^{p}\right).
\end{align*}
Applying the $L^{p}$ norm, we obtain that
\begin{equation}
\label{3.4}
\left\Vert Fy\right\Vert_{E^{\gamma,p}}\leq G_{\gamma}^{\ast}\left(
\left\Vert w\right\Vert _{1}+c\left(  1+\left(
\frac{1}{\Gamma(\gamma+1)}\right) ^{p}\right)R^{p}\right).
\end{equation}
Also,
\begin{align*}
\left\vert ^{C}D_{0}^{\gamma}L\left(  x\right)  \left(  t\right)\right\vert
&\leq\frac{3}{\Gamma\left(  2-\gamma\right)  }\int_{0}^{1} \left\vert
g\left(  s,x\left(  s\right)  \right)  \right\vert ds\\
& \quad +\frac{1}{\Gamma\left(  2-\gamma\right)  }\int_{0}^{1}\left\vert
q_{1}\left(  x\left(  s\right)  \right)  \right\vert ds+\frac{1}{\Gamma\left(
2-\gamma\right)  }\int_{0}^{1}\left\vert q_{2}\left(  x\left(  s\right)
\right)  \right\vert ds\\
&  \leq \frac{3}{\Gamma\left(  2-\gamma\right)  }\int_{0}^{1}\varphi_{1}\left(
s\right)  \left\vert x\left(  s\right)  \right\vert ds+\frac{1}{\Gamma\left(
2-\gamma\right)  }\int_{0}^{1}k_{1}  \left\vert x\left(
s\right)  \right\vert ds\\
& \quad +\frac{1}{\Gamma\left(  2-\gamma\right)  }
\int_{0}^{1}k_{2}   \left\vert x\left(  s\right)  \right\vert ds.
\end{align*}
Applying again the $L^{p}$ norm, we obtain from Holder's inequality that
\begin{equation*}
\left\Vert L\left(  x\right)  \right\Vert_{E^{\gamma,p}}
\leq\frac{3}{\Gamma\left(  2-\gamma\right)  }\left(
\left\Vert \varphi_{1}\right\Vert_{q}\left\Vert x\right\Vert _{p}\right)
+\frac{k_1}{\Gamma\left(2-\gamma\right)  }  \left\Vert x\right\Vert _{p}\\
+\frac{k_2}{\Gamma\left(  2-\gamma\right)  } \left\Vert x\right\Vert _{p}.
\end{equation*}
In view of \eqref{2.2}, we obtain
\begin{equation*}
\left\Vert L\left(  x\right)  \right\Vert _{E^{\gamma,p}}
\leq\left[\frac{3\left\Vert \varphi_{1}\right\Vert _{q}}{
\Gamma\left(  2-\gamma\right)\Gamma\left(  1+\gamma\right)}
+\frac{k_1}{\Gamma\left(  2-\gamma\right)
\Gamma\left(  1+\gamma\right)}
+\frac{k_2}{\Gamma\left(
2-\gamma\right)  \Gamma\left(  1+\gamma\right)  }\right]
\left\Vert x\right\Vert _{E^{\gamma,p}}\text{.}
\end{equation*}
Then,
\begin{equation}
\label{3.5}
\left\Vert L\left(  x\right)  \right\Vert _{E^{\gamma,p}}\leq
\frac{R\left[  3\left\Vert \varphi_{1}\right\Vert_{q}+k_1+k_2\right] }{
\Gamma\left(  2-\gamma\right)  \Gamma\left(  1+\gamma\right)  }.
\end{equation}
From \eqref{3.2}, \eqref{3.4} and \eqref{3.5},
we conclude that $Fx+Ly\in B_{R}$ whenever $x,y\in B_{R}$. \newline

\item[(iii)] Let us prove that $F\left(  B_{R}\right)
=\left\{  F\left(  u\right)  :u\in B_{R}\right\}$
is relatively compact in $E^{\gamma,p}$. Let $t\in\left(  0,1\right)$
and $h>0$, where $t+h\leq1$, and let $u\in D_{R}$. From \eqref{3.1}, we obtain that
\begin{align*}
&  \left\vert ^{C}D_{0}^{\gamma}Fy\left(  t+h\right)  -^{C}D_{0}^{\gamma}
Fy\left(  t\right)  \right\vert \\
&  \quad \leq\int_{0}^{1}\left\vert G_{\gamma}\left(  t+h,s\right)
-G_{\gamma}\left(  t,s\right)  \right\vert \left\vert f\left(  s,y\left(  s\right),
D^{\gamma}y\left(  s\right)  \right)  \right\vert ds\\
&  \quad \leq\int_{0}^{1}\left\vert G_{\gamma}\left(  t+h,s\right)
-G_{\gamma}\left(  t,s\right)  \right\vert \left[  w\left(  s\right)
+c\left(\left\vert y\left(  s\right)  \right\vert ^{p}+\left\vert D^{\gamma}
y\left(s\right)  \right\vert ^{p}\right)  \right]  ds\\
&  \quad \leq\sup_{t\in [0,1]}\left[  \sup_{s\in [0,1]}\left\vert G_{\gamma}\left(t+h,s\right)
-G_{\gamma}\left(  t,s\right)  \right\vert \right]  \left(\left\Vert
w\right\Vert _{1}+c\left(  1+\left( \frac{1}{\Gamma(\gamma+1)}\right) ^{p}\right)R^{p}\right).
\end{align*}
Therefore,
\begin{equation}
\label{3.10}
\frac{\left\Vert Fu\left(\cdot+h\right)
-Fu\left(\cdot \right)  \right\Vert_{E^{\gamma,p}}}{\left(
\left\Vert w\right\Vert _{1}+c\left(  1
+\left( \frac{1}{\Gamma(\gamma+1)}\right) ^{p}\right)R^{p}\right)}
\leq \sup_{t\in [0,1]}\left[  \sup_{s\in [0,1]}\left\vert G_{\gamma}\left(
t+h,s\right)  -G_{\gamma}\left(  t,s\right)  \right\vert \right].
\end{equation}
Then, $\left\Vert Fu\left(\cdot +h\right)-Fu\left(\cdot\right)
\right\Vert_{E^{\gamma,p}}\rightarrow0$ as $h\rightarrow 0$
for any $u\in B_{R}$, since $G_{\gamma}$ is a continuous function on
$\left[  0,1\right]  \times\left[0,1\right] $. From Lemma~\ref{lem3},
we conclude that $F:B_{R} \rightarrow B_{R}$ is compact. \newline

\item[(iv)] Finally, we prove that $L$ is a contraction. Let $x,y\in D_{R}$ and
$t\in\left(  0,1\right)$. Then,
\begin{align*}
\left\vert ^{C}D_{0}^{\gamma}L\left(  x\right)  \left(  t\right)
-^{C}D_{0}^{\gamma}L\left(  y\right)  \left(  t\right)  \right\vert
&  \leq\frac{3}{\Gamma\left(  2-\gamma\right)  }
\int_{0}^{1} \left\vert g\left(  s,x\left(  s\right)  \right)
-g\left(s,y\left(  s\right)  \right)  \right\vert ds\\
&  \quad +\frac{1}{\Gamma\left(  2-\gamma\right)  }\int_{0}^{1}\left\vert
q_{1}\left(  x\left(  s\right)  \right)  -q_{1}\left(  x\left(  s\right)
\right)  \right\vert ds\\
&  \quad +\frac{1}{\Gamma\left(  2-\gamma\right)  }\int_{0}^{1}\left\vert
q_{2}\left(  x\left(  s\right)  \right)  -q_{2}\left(  x\left(  s\right)
\right)  \right\vert ds\\
&  \leq\frac{3}{\Gamma\left(  2-\gamma\right)  }\int_{0}^{1}\varphi_{1}\left(
s\right)  \left\vert x\left(  s\right)  -y\left(  s\right)  \right\vert ds\\
&  \quad +\frac{k_1}{\Gamma\left(  2-\gamma\right)  }
\int_{0}^{1} \left\vert x\left(  s\right)  -x\left(  s\right)  \right\vert ds\\
&  \quad +\frac{k_2}{\Gamma\left(  2-\gamma\right)  }\int_{0}^{1}
\left\vert x\left(  s\right)  -x\left(  s\right)  \right\vert ds.
\end{align*}
Applying the $L^{p}$ norm and Holder's inequality, we obtain that
\begin{gather*}
\left\Vert L\left(  x\right)  -L\left(  y\right)  \right\Vert_{E^{\gamma,p}}
\leq\frac{3}{\Gamma\left(  2-\gamma\right)  }\left(
\left\Vert \varphi_{1}\right\Vert _{q}\left\Vert x-y\right\Vert _{p}\right)
+\frac{k_1}{\Gamma\left(  2-\gamma\right)  }\left(\left\Vert x-y\right\Vert _{p}\right)\\
+\frac{k_2}{\Gamma\left(  2-\gamma\right)  }\left(
\left\Vert x-y\right\Vert_{p}\right).
\end{gather*}
Then, from \eqref{2.2}, we obtain
\[
\left\Vert L\left(  x\right)  -L\left(  y\right)
\right\Vert_{E^{\gamma,p}}
\leq\frac{3\left\Vert \varphi_{1}\right\Vert_{q}
+k_1+k_2}{\Gamma\left(2-\gamma\right)  \Gamma\left(
1+\gamma\right)}\left\Vert x-y\right\Vert_{E^{\gamma,p}}.
\]
From \eqref{3.3}, the operator $L$ is a contraction. 
\end{enumerate}
As a consequence of (i)--(iv), we conclude that $F:B_{R}\rightarrow B_{R}$ 
is continuous and compact.  As a consequence of Krasnosel'skii fixed point 
theorem, we deduce that $F+G$ has a fixed point $y\in B_{R}\subset E^{\gamma,p}$,
which is a solution to problem \eqref{77}.
\end{proof}

We now illustrate our Theorem~\ref{the3.1} with {two examples}.

\begin{Example}
Consider the fractional boundary value problem \eqref{77} with
\begin{align*}
\alpha &=2.5, \quad \gamma=0.5, \quad p=3, \quad q=\frac{3}{2},\\
f\left(  t,x,y\right)   &  =\frac{\exp\left(  -t\right)  }{5}-\frac{1}{164\pi}
\arctan\left(  x^{3}+y^{3}\right),\\
g\left(  t,x\right) &=\frac{1}{10}t^{\frac{2}{3}}x,\\
q_{1}\left(  x\right) &=q_{2}\left( x\right)  =\frac{1}{20}x,
\end{align*}
which we denote by {$(P_1)$}. Hypotheses $\left(  H1\right)$
and $\left(  H2\right)$ are satisfied for 
$$
{w\left(  t\right)
=\frac{\exp\left(  -t\right)  }{5}\in L^{1}\left(  0,1\right),
\quad c=\frac{1}{164\pi}, \quad 
\varphi_{1}\left(  t\right) =\frac{t^{\frac{2}{3}}}{10}, 
\quad \text{and} \quad k_1=k_2 =\frac{1}{20}.}
$$
Moreover, we have
\[
\frac{\left[  3\left\Vert \varphi_{1}\right\Vert _{q}+k_1+k_2\right]}
{\Gamma\left(  2-\gamma\right)  \Gamma\left(  1+\gamma\right)}
=\frac{\left[
\frac{3}{10}\left( \frac{1}{2}\right)^{\frac{2}{3}}
+\frac{1}{10}\right]  }{\left( \Gamma\left(
\frac{3}{2}\right)\right) ^{2}}
\simeq 0.368<1\text{.}
\]
If we choose $R=2$, then we obtain
\begin{align*}
&  \frac{R\left[  3\left\Vert \varphi_{1}\right\Vert _{q}+k_1+k_2\right]}{
\Gamma\left(  2-\gamma\right)  \Gamma\left(  1+\gamma\right)  }+G_{\gamma}^{\ast}\left(
\left\Vert w\right\Vert _{1}+c\left(  1+\left( \frac{1}{\Gamma(\gamma+1)}\right) ^{p}\right)
R^{p}\right) -R\\
& \hspace{0.3cm} \leq \frac{2\left[  \frac{3}{10}\left(\frac{1}{2} \right)^{
\frac{2}{3}} +\frac{1}{10}\right]  }{\left( \Gamma\left(
\frac{3}{2}\right) \right)^{2}} +  4.047 \left(  \frac{1}{5}
+\frac{1}{164\pi}\left(  1+\left( \frac{1}{\Gamma\left(
\frac{3}{2}\right)  }\right) ^{3}\right)  2^{3}\right) -2\\
& \hspace{0.3cm} \simeq -0.301.
\end{align*}
Since all conditions of our Theorem~\ref{the3.1} are satisfied,
we conclude that the fractional boundary value problem {$(P_1)$}
has a solution in $E^{\gamma,p}$.
\end{Example}

\begin{Example}
Consider the fractional boundary value problem \eqref{77} with
\begin{align*}
\alpha &=2.7, \quad \gamma=0.7, \quad p=4, \quad q=\frac{4}{3},\\
f\left(  t,x,y\right)   &  =\frac{1}{10}\sin\left(  t\right) + \frac{1}{200}
\cos\left(  x^{4}+y^{4}\right),\\
g\left(  t,x\right) &=\frac{1}{9\pi}t^{\frac{3}{4}}\arctan(x),\\
q_{1}\left(  x\right) &=q_{2}\left( x\right)  =\frac{1}{10}\sin(x),
\end{align*}
which we denote by  {$(P_2)$}. Hypotheses $\left(  H1\right)$
and $\left(  H2\right)$ are satisfied for 
$$
w\left(  t\right)=\frac{1}{10}\sin\left( t\right)\in L^{1}\left(  0,1\right),
\quad c=\frac{1}{200}, \quad 
\varphi_{1}\left(  t\right)=\frac{t^{\frac{3}{4}}}{9\pi} 
\quad \text{and} \quad k_1=k_2 =\frac{1}{10}. 
$$
Moreover, we have
\[
\frac{\left[  3\left\Vert \varphi_{1}\right\Vert _{q}+k_1+k_2\right]}
{\Gamma\left(  2-\gamma\right)  \Gamma\left(  1+\gamma\right)}
=\frac{\left[
\frac{1}{3\pi}\left( \frac{1}{2}\right)^{\frac{3}{4}}
+\frac{1}{5}\right]  }{ \Gamma(1.3)\Gamma(1.7) }
\simeq 0.323<1\text{.}
\]
If we choose $R=2$, then we obtain
\begin{align*}
& \frac{R\left[  
3\left\Vert \varphi_{1}\right\Vert _{q}+k_1+k_2\right]}{
\Gamma\left(  2-\gamma\right)  \Gamma\left(  1+\gamma\right)  }+G_{\gamma}^{\ast}\left(
\left\Vert w\right\Vert _{1}+c\left(  1+\left( \frac{1}{\Gamma(\gamma+1)}\right) ^{p}\right)
R^{p}\right) -R\\
& \hspace{0.3cm} \leq \frac{2\left[  \frac{1}{3\pi}\left(\frac{1}{2} \right)^{
\frac{3}{4}} +\frac{1}{5}\right]  }{\Gamma(1.3)\Gamma(1.7)} +  3.9995 \left(  \frac{1}{10}
+\frac{1}{200}\left(  1+\left( \frac{1}{\Gamma\left(  1.7 \right)  }\right) ^{4}\right)  2^{4}\right) -2\\
& \hspace{0.3cm} \simeq -0.163.
\end{align*}
Since all conditions of our Theorem~\ref{the3.1} are satisfied,
we conclude that the fractional boundary value problem {$(P_2)$}
has a solution in $E^{\gamma,p}$.
\end{Example}


\section{Discussion}
\label{sec:disc}

The celebrated existence result of Ahmad and Nieto \cite{14}
for problem \eqref{66}  {is} obtained via Krasnosel'skii fixed-point
theorem in the space of continuous functions. For that, they needed
to apply Ascoli's theorem in order to provide the compactness
of the first part of the Krasnosel'skii operator. Here, we proved
existence for the more general problem \eqref{77} in the
fractional derivative Banach space $E^{\gamma,p}$,
equipped with the norm \eqref{n}. From norm \eqref{n}, it is natural to
deal with a subspace of $L^{p}\times L^{p}$. Since Ascoli's theorem is limited
to the space of continuous functions for the compactness, we had to make use
of a different approach to ensure existence of solution in the fractional derivative
space $E^{\gamma,p}$. Our tool was the Kolmogorov--Riesz compactness theorem,
which turned out to be a powerful tool to address the problem. To the best of our
knowledge, the use of the Kolmogorov--Riesz compactness theorem to prove
existence results for boundary value problems involving nonlinear integrodifferential
equations of fractional order with integral boundary conditions is a completely new approach.
In this direction, we are only aware of the work \cite{MR3856969}, where a necessary
and sufficient condition of pre-compactness in variable exponent Lebesgue spaces is established
and, as an application, the existence of solutions to a fractional Cauchy problem is obtained
in the Lebesgue space of variable exponent. As future work, we intend to generalize
our existence result to the variable-order case \cite{MR3822307}.


\vspace{6pt}


\authorcontributions{Conceptualization, D.B., A.C. and D.F.M.T.; 
validation, D.B., A.C. and D.F.M.T.; 
formal analysis, D.B., A.C. and D.F.M.T.; 
investigation, D.B., A.C. and D.F.M.T.; 
writing---original draft preparation, D. B., A.C. and D.F.M.T.; 
writing---review and editing, D.B., A.C. and D.F.M.T.
All authors have read and agreed to the published version of the manuscript.}
	
\funding{This research was partially funded by FCT, grant number UIDB/04106/2020 (CIDMA).}

\institutionalreview{Not applicable.}

\informedconsent{Not applicable.}

\dataavailability{Data sharing not applicable.}

\acknowledgments{{The authors are grateful to the referees for their comments and remarks.}}
	
\conflictsofinterest{The authors declare no conflict of interest.}


\end{paracol}

\reftitle{References}


\end{document}